\documentclass[12pt]{article}
\usepackage{amsmath,amsfonts,amscd,amssymb,latexsym}

\title{Algebraic Geometry over Free Metabelian  Lie Algebra I: $U$-Algebras and Universal Classes}
\author{\textsf{E. Yu. Daniyarova}\thanks{The first author was supported by
the RFFI grant N02-01-00192.} \and \textsf{I. V.
Kazachkov}\thanks{The first and the second authors were supported
by the `Universitety Rossii' grant.} \and \textsf{V. N.
Remeslennikov}$^{\dag}$}

\newtheorem{lem}{Lemma}[section]
\newtheorem{thm}{Theorem}[section]
\newtheorem{cor}{Corollary}[section]
\newtheorem{prop}{Proposition}[section]

\newtheorem{defn}{Definition}[section]

\newtheorem{rem}{Remark}[section]

\newcommand{\Beta}{\mathcal{B}}

\newcommand{\F}{F}

\newcommand{\proof}{\paragraph{Proof.}}
\renewcommand{\hom}{\texttt{Hom}}

\newcommand{\qvar}{\texttt{qvar}}
\newcommand{\var}{\texttt{var}}

\newcommand{\Fit}{\texttt{Fit}}
\newcommand{\ucl}{\texttt{ucl}}

\begin{document}

\maketitle

\begin{abstract}
This paper is the first in a series of three, the aim of which
is to lay the foundations of algebraic geometry over the free
metabelian Lie algebra $F$. In the current paper we introduce the
notion of a metabelian Lie $U$-algebra and establish connections
between metabelian Lie $U$-algebras and special matrix Lie
algebras. We define the $\Delta $-localisation of a metabelian
Lie $U$-algebra $A$ and the direct module extension of the
Fitting's radical of $A$ and show that these algebras lie in the
universal closure of $A$.
\end{abstract}

\tableofcontents

\section{Introduction}
This paper is the first in a series of three, the object of which
is to lay the foundations of algebraic geometry over free
metabelian Lie algebra $F$. In papers \cite{alg1}, \cite{alg2} G.
Baumslag, A. Myasnikov and V. Remeslennikov  introduced and
studied group-theoretic counterparts to the main notions of
classical algebraic geometry. The results of \cite{alg1},
\cite{alg2} can be extended to an arbitrary algebraic system.

The main aim of the current series of our papers is to create a
structural theory of Lie algebras from the quasivariety $\qvar(F)$
and the universal closure $\ucl(F)$ generated by the free
metabelian Lie algebra $F$. This motivation is due to the
following circumstance: finitely generated Lie algebras from the
quasivariety $\qvar(F)$ (the universal closure $\ucl(F)$) are
exactly the coordinate algebras of (irreducible) algebraic sets
over $F$.

This paper holds some preliminary results, which will play an
important role for algebraic geometry over $F$. In Section
\ref{sec:2} we list some preliminary results on metabelian Lie
algebras.  Some of these results are known, but we give new proofs
and further use not only the statements itself, but the methods of
proofs.

In papers of V. A. Artamonov \cite{Art1} and \cite{Art2} the
author gives a presentation of the free metabelian Lie algebra in
a module over the ring of polynomials. Note that this
presentation is a particular case of a more general construction
of embedding a Lie algebras of special type into verbal wreath
product of Lie algebras. The latter construction is due to A. L.
Shmel'kin (see \cite{Shm}). In Section \ref{sec:3} we borough the
ideas of the authors mentioned to construct special matrix
metabelian Lie algebras. There we introduce the notion of a
$U$-algebra, which is very important for the use of algebraic
geometry over the free metabelian Lie algebra $F$ and establish
their connections with special matrix Lie algebras. In particular
(see Theorem \ref{thm:321}), every finitely generated $U$-algebra
is a subalgebra of a special matrix Lie algebra.

In Section \ref{sec:4} for a fixed Lie algebra $A$ we introduce
the first order language  $L_A $ and study several universal
classes in this language (the universal closure, the quasivariety,
etc.). In this section we also introduce the notions of the
$\Delta $-localisation of a $U$-algebra $A$ and of the direct
module extension of the Fitting's radical of the algebra $A$. We
show (see Propositions \ref{prop:422} and \ref{prop:435}) that the
new algebras lie in the universal closure of the algebra $A$. We
refer to \cite{Bah} for  preliminaries of Lie algebras and to
\cite{Bourb} and \cite{Leng} for commutative algebra.

Some of the results of this paper can be found in author's
preprint \cite{preprint}.

\section{Metabelian Lie Algebras} \label{sec:2}
In this section we give an exposition of preliminary results,
which will be used further. Some of these results are well-known.

Recall that a vector space $A$ over a field  $k$ equipped with a
bilinear multiplication satisfying  the following universal axioms
\begin{itemize}
    \item  The anti-commutativity identity $a \circ a = 0$
    (notice that the regular anti-commutativity
    axiom $a \circ b =  - b \circ a$ is implied by  $a \circ a =
    0$),
    \item The Jacoby identity $(a \circ b) \circ c + (b \circ c) \circ a + (c \circ a)
            \circ b = 0$,
\end{itemize}
is termed a Lie algebra.

In the event that the Lie algebra $A$ satisfies
\begin{itemize}
    \item The metabelian identity $(a \circ b) \circ (c \circ d) =
    0$;
\end{itemize}
it is termed a metabelian Lie algebra.

The regular notation of multiplication of elements $a$ and $b$
from $A$ is `$\circ$', $a \circ b$. For the sake of brevity, we
sometimes omit `$\circ$' and use the notation $ab$.

Left-normed products $a_1 a_2 a_3 \cdots a_n$, of $a_1 ,a_2 ,a_3
,\ldots,a_n \in A$ are defined as
$$
(\ldots((a_1  \circ a_2 ) \circ a_3 ) \circ \ldots) \circ a_n.
$$
We term such products by left-normed words or monomials of the
degree or of the length  $n$. Similar, right-hand side definition
of such words gives rise to the notion of a right-normed monomial.
As it is well-known in the theory of Lie algebras every element of
degree $n$ from letters $a_1 ,\ldots, a_n $ can be presented as a
linear combination of left-normed (right-normed) monomials of
length $n$ from these letters.

\subsection{The Fitting's Radical and the Commutant}

Two special ideals of a metabelian Lie algebra $A$ are of great
importance for us.

\begin{defn}[The commutant]
The ideal generated by the set $\left\{ a \circ b| \; a,b \in A
\right\}$ is termed the commutant of the algebra  $A$ and denoted
by $A^2 $. In the event that $A$ is a metabelian Lie algebra,
resulting from the metabelian identity, the commutant $A^2 $ is
abelian.
\end{defn}

\begin{defn}[The Fitting's radical] \label{defn:II-fit}
The ideal generated by the elements from nilpotent ideals of the
algebra $A$  is termed the Fitting's radical and denoted by {\rm
$\Fit (A)$}.
\end{defn}

Since the commutant of a metabelian Lie algebra is an Abelian
ideal, it is always contained in the Fitting's radical, $A^2
\subseteq \Fit(A)$.

\begin{rem}
For this notation, there exists a metabelian Lie algebra $A$ such
that {\rm $A^2 \ne \Fit(A)$}. If $A$ is Abelian then {\rm $A =
\Fit(A)$} while $A^2  = 0$. If $A \ne 0$, then \/ {\rm $\Fit(A)
\ne 0$}.
\end{rem}

Metabelian Lie algebras have the following useful property. Every
left-normed monomial $abc_1 \cdots c_n $, where here $n \ge 2$
equals $ab c_{\tau (1)} \cdots c_{\tau (n)}$, where $\tau$ is a
permutation on the set of indices $\left\{1, \ldots, n\right\}$.

To prove that a permutation of the far co-factors preserves the
element, it suffices to show that a permutation of two adjacent
co-factors $c_i$ and $c_{i + 1}$ preserves the element. Consider
the following equation
$$
abc_1 \cdots c_i c_{i + 1} \cdots c_k  - abc_1 \cdots c_{i + 1}
c_i \cdots c_k = (abc_1 \cdots c_i c_{i + 1}  - abc_1 \cdots c_{i
+ 1} c_i )\cdots c_k .
$$
Using the anti-commutativity identity: $abc_1 \cdots c_{i + 1} c_i
= - c_{i + 1} (abc_1 \cdots )c_i $. By the metabelian identity:
$c_i c_{i + 1} (abc_1 \ldots) = 0$. The initial expression,
therefore, rewrites as follows
$$
((abc_1 \cdots )c_i c_{i + 1}  + c_{i + 1} (abc_1 \cdots )c_i  +
c_i c_{i + 1} (abc_1 \cdots))\cdots c_k.
$$
The latter equals zero by the Jacoby identity.

We shall make use of some facts on nilpotent subalgebras of
metabelian Lie algebras.

\begin{lem} \label{lem:nilpid}
Let $A$ be a metabelian Lie algebra. And let $I_1$ and $I_2 $ be
nilpotent ideals of $A$ of nilpotency classes $n_1$ and $n_2 $.
Then the ideal $I =  \left< I_1 ,I_2  \right> $ is also nilpotent.
\end{lem}
\proof It is fairly easy to see that the nilpotency class of $I$
is lower than or equals $2n$, where $n = \max \left\{ n_1 ,n_2
\right\}$. \hfill $\blacksquare$

\begin{cor}
Let $I_1,I_2, \ldots, I_n $ be nilpotent ideals of $A$. Then the
ideal $I = \left< I_1 ,I_2, \ldots I_n  \right> $ is nilpotent.
\end{cor}

\begin{cor} \label{cor:II-fitchar}
An element $x \in A$ lies in the Fitting's radical \/ {\rm
$\Fit(A)$} if and only if $I = \left< x \right> $ is nilpotent.
Or, which is equivalent, if  $x$ lies in some nilpotent ideal of
$A$.
\end{cor}

\begin{lem}
Let $C$ be a nilpotent metabelian Lie algebra of nilpotency class
$n \ge 2$. Then  $C$ contains a 2-generated nilpotent subalgebra
$D$ of nilpotency class 2.
\end{lem}
\proof Let
$$
\left\{ 0 \right \}  = Z_0 (C) < Z_1 (C) < Z_2 (C) < \ldots < Z_{n
- 1} (C) < Z_n (C) = C
$$
be the upper central series for $C$ (see \cite{Bah}). Since $C$ is
non-Abelian, the set $Z_2 (C)\smallsetminus Z_1 (C)$ is non-empty.
Take $c_1  \in Z_2 (C)\smallsetminus Z_1 (C)$ and $c_2  \in C$ so
that $c_1  \circ c_2  \ne 0$. The algebra $D =  \left< c_1 ,c_2
\right>$ is a 2-generated nilpotent Lie algebra of class 2. \hfill
$\blacksquare$

\begin{lem}
If every nilpotent ideal of $A$ is Abelian then \/ {\rm $\Fit(A)$}
is Abelian.
\end{lem}
\proof Assume the converse. There, therefore, exist two elements
$c_1,c_2 \in A$ and two nilpotent ideals $I_1$ and $I_2 $ so that
$c_1 \in I_1$ and $c_2 \in I_2$  and so that $c_1  \circ c_2  \ne
0$. According to Lemma \ref{lem:nilpid}, the ideal $I =  \left<
I_1 ,I_2  \right> $ is nilpotent but not Abelian. This derives a
contradiction. \hfill $\blacksquare$

\begin{lem} \label{lem:214}
Let $A$ be a metabelian Lie algebra and suppose that the element
$a \in A$ commutes with every element from $A^2 $. Then {\rm $a
\in \Fit(A)$}.
\end{lem}
\proof Straightforward checking shows that under the assumptions
of the lemma the principal ideal $I = \left< a \right>$ is
Abelian. Consequently, by Corollary \ref{cor:II-fitchar}, $a \in
\Fit(A)$. \hfill $\blacksquare$

This series of lemmas allows to introduce a structure of a module
over the ring of polynomials on the commutant $A^2$ and on the
Fitting's radical $\Fit (A)$ (in the event that it is Abelian).
Or, more generally, a structure of a module on an arbitrary
abelian ideal $I$ of $A$, containing $A^2$.

In that case ${\raise0.7ex\hbox{$A$} \!\mathord{\left/
 {\vphantom {A I}}\right.\kern-\nulldelimiterspace}
\!\lower0.7ex\hbox{$I$}}$ is also Abelian, since $A^2 \subseteq
I$. Take a maximal linearly independent  modulo $I$ system of
elements $\left\{ a_\alpha ,\alpha \in \Lambda \right\} $ from
$A$. The ideal $I$ admits a structure of a module over the ring of
polynomials $R = k\left[x_\alpha ,\alpha \in \Lambda \right]$. The
addition and the multiplication by elements of $k$ is correctly
defined, since $I$ is a vector space over $k$. By the definition,
for an element $b \in I$,  set
$$ b \cdot x_\alpha   = b \circ
a_\alpha ,\;\alpha  \in \Lambda.
$$

 To define the action of a singleton $f$ from  $R$ on the
element $b \in I$ we use induction on the degree $\deg f$. By the
definition, $b \cdot (fx_\alpha) = (b \cdot f) \cdot x_\alpha$.
Let $f_1, \ldots, f_m$ be a tuple of singletons from $R$ and let
$\gamma_1,\ldots, \gamma_m \in k$. We set
\begin{equation} \label{eq:II-act} b \cdot
(\gamma_1 f_1+\cdots+\gamma_mf_m)=\gamma_1 b \cdot
f+\cdots+\gamma_m b \cdot f_m.
\end{equation}

Equation (\ref{eq:II-act}) defines a structure of a module on the
ideal $I$ over the ring of polynomials with non-commuting
variables. We next show that this definition is correct for the
ring of commutative polynomials $R$. It is sufficient to show that
$b \cdot (x_\alpha x_\beta ) = b \cdot (x_\beta  x_\alpha )$,
where here $b \in I$ and $\alpha ,\beta \in \Lambda $ is an
arbitrary pair of indices. Applying the Jacoby identity we obtain
$$
b \cdot (x_\alpha  x_\beta  ) - b \cdot (x_\beta  x_\alpha  )
    = b \circ a_\alpha   \circ a_\beta- b \circ a_\beta   \circ a_\alpha =
    a_\beta   \circ a_\alpha   \circ b.
$$
Since ${\raise0.7ex\hbox{$A$} \!\mathord{\left/
 {\vphantom {A I}}\right.\kern-\nulldelimiterspace}
\!\lower0.7ex\hbox{$I$}}$ and $I$ are Abelian, $a_\beta \circ
a_\alpha \in I$  and  $a_\beta   \circ a_\alpha   \circ b = 0$,
from the foregoing argument, we obtain that $I$ is a module over
the ring $R$.

\begin{rem}
How does the module structure of $I$ change if instead of the
initial maximal linearly independent modulo $I$ system $\left\{
a_\alpha ,\alpha \in \Lambda \right\}$ of elements from $A$ we
choose another such system $\left\{ a'_\alpha ,\alpha \in \Lambda
\right\}$? Since the ideal $I$ is Abelian such transformation
implies only a $k$-linear change of the variables from the ring of
polynomials $R$.
\end{rem}

\begin{rem} \label{rem:dirsum}
Denote by  $V$ the $k$-linear span of the set
 $\left\{ a_\alpha  ,\alpha  \in \Lambda \right\}$. Then the algebra
 $A$, treated as a vector space, decomposes into the following direct sum $V \oplus I$
of vector spaces over $k$.
\end{rem}

\begin{rem} \label{rem:II-endom}
Furtherafter we shall sometimes  abuse the notation $b \cdot f$,
where $f \in R$ and $b \in A$ (not necessarily $b \in I$). By this
notation we mean that if
$$
f = \gamma  + x_{\alpha _1 } f_1  + \cdots + x_{\alpha _l } f_l,
\quad \gamma  \in k, \; f_i  \in R, \; \alpha  \in \Lambda
$$
then
\begin{equation} \label{eq:II:f}
b \cdot f = \gamma b + (ba_{\alpha _1 } ) \cdot f_1  + \cdots +
(ba_{\alpha _l } ) \cdot f_l.
\end{equation}
This equation is correctly  defined, since $ ba_{\alpha _i } \in
I$. Notice that Equation {\rm (\ref{eq:II:f})} depends on the
choice of designated letters $x_{\alpha _1 }, \ldots ,x_{\alpha
_l}$. We use this notation in the context when this choice is not
significant.
\end{rem}

\subsection{Generators and Defining Relations of Metabelian Lie
Algebras and Finitely Generated Metabelian Lie Algebras} \label{sec:21}

Below we describe a convenient set of generators of an arbitrary
metabelian Lie algebra $A$.

As above, let $A$ be a metabelian Lie algebra and let $I$ be an
ideal of $A$ such that ${\raise0.7ex\hbox{$A$} \!\mathord{\left/
 {\vphantom {A I}}\right.\kern-\nulldelimiterspace}
\!\lower0.7ex\hbox{$I$}}$ is Abelian. We usually think of $I$ as
 the commutant or the Fitting's radical (in the event that it is Abelian).
 Take maximal linearly independent modulo $I$
system  of elements $\left\{ a_\alpha ,\alpha \in \Lambda
\right\}$ from $A$. Denote by $V$ the $k$-linear span of $\left\{
a_\alpha ,\alpha \in \Lambda \right\}$. Suppose that $\left\{
b_\beta ,\beta  \in \Beta \right\}$ generates $I$ as an $R$-module
and is minimal.

\begin{lem} \label{lem:cansetgen}
The union of the sets $\left\{ a_\alpha ,\alpha \in \Lambda
\right\}$ and $\left\{ b_\beta ,\beta  \in \Beta \right\}$
generates the algebra $A$.
\end{lem}
\proof On account of Remark \ref{rem:dirsum}, every element $a \in
A$ can be written as $a = c + b$, where $c \in V$ and $b \in I$.
The element $b \in I$ has a presentation:
$$
b = b_{\beta _1 } \cdot f_1  + \cdots + b_{\beta _l }  \cdot f_l,
$$
where here $f_i  \in R$. Every summand $b_{\beta _i } \cdot f_i $,
treated as an element of $A$, is a Lie polynomial from $\left\{
a_\alpha, \alpha \in \Lambda \right\},\;b_{\beta _i }$.
Consequently,
 $b$ and, therefore, $a$  are Lie polynomials from
 $\left\{ a_\alpha ,\alpha \in \Lambda \right\} \cup \left\{ b_\beta
,\beta  \in \Beta \right\}$. \hfill $\blacksquare$

\begin{rem}
The system of generators from Lemma { \rm \ref{lem:cansetgen} } is
termed the canonical system of generators of a metabelian Lie
algebra. We shall use it in the following, when working with
metabelian Lie algebras.
\end{rem}

\begin{lem} \label{lem:222}
If  $A$ is a finitely generated metabelian Lie algebra then its
canonical set of generators is finite.
\end{lem}
\proof Let $\left\{ s_1 ,\ldots ,s_n \right\}$ be a finite set of
generators of the algebra $A$. We show that the canonical set of
generators is also finite.

First of all observe that ${\raise0.7ex\hbox{$A$}
\!\mathord{\left/
 {\vphantom {A I}}\right.\kern-\nulldelimiterspace}
\!\lower0.7ex\hbox{$I$}}$ is a finite dimensional vector space
over $k$. This is immediate, since every element of $A$ is a Lie
polynomial from the letters $\left\{ s_1 ,\ldots,s_n \right\}$.
Every polynomial admits the following presentation
$$
a = \alpha _1 s_1  + \cdots + \alpha _n s_n  + b.
$$
Here $\alpha _i  \in k$, and $b \in A^2$ and therefore $b \in I$.
In what follows that $\dim ({\raise0.7ex\hbox{$A$}
\!\mathord{\left/
 {\vphantom {A I}}\right.\kern-\nulldelimiterspace}
\!\lower0.7ex\hbox{$I$}}) = r$ is finite and, moreover, is lower
than $n$.

We next show that the ideal $I$, regarded as a module over the
ring $R = k\left[x_1 ,\ldots ,x_r \right]$ is finitely generated.
The commutant $A^2 $, treated as a $k$-vector space, is generated
by left-normed monomials $s_{i_1 }  \circ s_{i_2 }  \circ s_{i_3 }
\circ \cdots \circ s_{i_m }$ of the length $\ge 2$. Write $s_i $,
$i = 1,\ldots,n$, as a sum: $s_i  = c_i  + b_i $ , $c_i \in V$ ,
$b_i \in I$. Substituting and evaluating we see that $s_{i_1 }
\circ s_{i_2 }  \circ \cdots \circ s_{i_m } = (s_{i_1 } \circ
s_{i_2 } ) \cdot f$, where $f \in R$. So $A^2 $ as a module over
$R$ is generated by its finite subset $\left\{  {s_i  \circ s_j }
\mid \;i < j = 1,\ldots, n \right\}$. Since
${\raise0.7ex\hbox{$I$} \!\mathord{\left/
 {\vphantom {I {A^2 }}}\right.\kern-\nulldelimiterspace}
\!\lower0.7ex\hbox{${A^2 }$}}$ is also a finite dimensional
$k$-vector space, the ideal $I$ is a finitely generated
$R$-module.

We thereby proved that both parts of the canonical system of
generators are finite. \hfill $\blacksquare$

In the category of all metabelian Lie algebras over the field  $k$
every finitely generated $k$-algebra $A$ is finitely presented.
The proof of the next theorem establishes the form of the defining
relations of a finitely generated metabelian Lie algebra $A$.

\begin{thm} \label{thm:223}
Every finitely generated metabelian Lie  $k$-algebra $A$ is
finitely presented in the category of all metabelian Lie
$k$-algebras.
\end{thm}
\proof Consider the vector space ${\raise0.7ex\hbox{$A$}
\!\mathord{\left/ {\vphantom {A {A^2
}}}\right.\kern-\nulldelimiterspace} \!\lower0.7ex\hbox{${A^2
}$}}$. Suppose that  $\dim {\raise0.7ex\hbox{$A$}
\!\mathord{\left/ {\vphantom {A {A^2
}}}\right.\kern-\nulldelimiterspace} \!\lower0.7ex\hbox{${A^2
}$}}=r$. The commutant $A^2$ is a finitely generated module over
the ring $R = k\left[x_1 ,\ldots,x_r \right]$. Let
$b_1,\ldots,b_l$ be its module generators. We shall use the
canonical system of generators of $A$: $a_1 ,\ldots,a_r ,\;b_1
,\ldots,b_l$.

Since the ring $R$ is Noetherian and since $A^2 $ is a finitely
generated $R$-module, the module $A^2 $ is finitely presented (see
\cite{Bourb, Leng}). Fix a finite presentation of the module $A^2
$. An arbitrary relation of the module $A^2$ has the form:
$$
b_1  \cdot f_1  + \cdots + b_l  \cdot f_l  = 0,\;f_i  \in R.
$$
This relation rewrites into the signature of metabelian Lie
algebras. This way, we obtain a relation, which involves letters
$a_1,\ldots,a_r$, $b_1 ,\ldots,b_l $ and which is true in the
algebra $A$. By the definition, the relations of this kind form
the first type of defining relations of $A$.

There are three types of relations in $A$. The second type are the
obvious relations, which show that the commutant is Abelian:
$$
b_i\circ b_j  = 0,\;\,b_i  \circ a_s  \circ b_j  = 0,\quad \,i < j
= 1,\ldots,l,\;\,s = 1,\ldots,r.
$$

The third type of the relations is to show that for every pair of
indices  $(i,j)$ the product  $a_i  \circ a_j $ lies in the
commutant:
$$
a_i  \circ a_j  = b_1  \cdot g_{ij} ^1  + \cdots + b_l \cdot
g_{ij} ^l ,\quad i < j = 1,\ldots,r,\;\,g_{ij} ^p  \in R.
$$
Again, these relations rewrite into the signature of metabelian
Lie algebras and then form the third type of defining relations
for $A$.

We are now left to show that the algebra $A$ is completely defined
by the relations of these three types. Let  $F$ be the free
metabelian Lie algebra with the base $x_1 ,\ldots,x_r$, $y_1
,\ldots,y_l $. And let $K$ be the ideal of $F$ generated by the
relations of $A$, which are one of the three types described
above. Denote by $\varphi :F \rightarrow A$ the canonical
epimorphism, defined by the following equalities: $\varphi (x_i )
= a_i ,\;i = 1,\ldots,r$, $\varphi (y_j ) = b_j ,\;j =
1,\ldots,l$. We next show that $ \ker \varphi  = K$. The inclusion
$\ker\varphi \supseteq K$ is obvious. Only the converse inclusion
is at issue.

Consider an arbitrary Lie polynomial $u$ from the letters $x_1
,\ldots,x_r$, $y_1 ,\ldots,y_l $. Since $u$ is an element of a Lie
$k$-algebra, it rewrites as a $k$-linear combination of
left-normed monomials from the same letters. As shown above, a
permutation of the $c_i$'s in left-normed monomials of the form
$abc_1 \cdots c_n$ preserves the element (here $a,b, c_1, \ldots,
c_n \in \left\{x_1 ,\ldots,x_r, y_1 ,\ldots,y_l\right\}$).

Using the  relations of the third type and the anti-commutativity
identity we write every monomial of the length $ \ge 2$ as a
linear combination of monomials, the initial letter of which is
one of the $y_j$'s. Now, if there are at least two occurrences of
the $y_j$'s in a monomial then, on account of the second type of
the relations, this monomial lies in $K$. In what follows that an
arbitrary element $u \in F$ can be presented in the following
form:
$$
u = \alpha _1 x_1  + \cdots + \alpha _r x_r  + y_1  \cdot f_1  +
\cdots + y_l  \cdot f_l  + v,
$$
where $v \in K$, $\alpha _i  \in k$, $f_j  \in R$. Suppose that $u
\in \ker \varphi$. We, therefore, obtain:
$$
\alpha _1 a_1  + \cdots + \alpha _r a_r  + b_1  \cdot f_1  +
\cdots + b_l \cdot f_l  = 0.
$$
Consequently, $\alpha _1  = \ldots = \alpha _r  = 0$ and  $b_1
\cdot f_1  + \cdots + b_l  \cdot f_l  = 0$. The latter equality
and the first type of the relations together imply that $y_1 \cdot
f_1 + \cdots + y_l \cdot f_l \in K$. And we, therefore, have shown
that $u \in K$. \hfill $\blacksquare$

\begin{rem} We suppose that the result of Theorem {\rm
\ref{thm:223}} is known, although we have not found any reference.
Furthermore, in the following we use the particular form of the
relations given in the proof of Theorem {\rm \ref{thm:223}}.
\end{rem}

\begin{cor} In the event that the Fitting's radical \/ {\rm $\Fit(A)$} is Abelian
the presentation given in the proof of Theorem {\rm \ref{thm:223}}
can be refined. The new presentation is obtained using the
argument of Theorem {\rm \ref{thm:223}} with $A^2 $ replaced by
\/{\rm $\Fit(A)$}.
\end{cor}

\begin{rem}
For infinitely generated metabelian Lie algebras the three types
of relations introduced in the proof of Theorem {\rm
\ref{thm:223}} are defining. Though, any of the three types may be
infinite.
\end{rem}

We next introduce a class of endomorphisms of a metabelian Lie
algebra which will play an important role in the following.
Suppose that $I$ is either the commutant  $A^2$ or  the Fitting's
radical $\Fit(A)$ (in the event that it is Abelian).

Let $f \in R$ be a polynomial with zero free term. We define the
map $\varphi$ as the multiplication of elements of the algebra $A$
by the polynomial $(1+f)$. Recall that  this map is correctly
defined only for the elements of $I$. To extend the action of
$\varphi$ to the algebra $A$ we write the polynomial $f$ in the
form $f = x_{j_1 } f_1  + \cdots + x_{j_q } f_q$ (see Remark
\ref{rem:II-endom}). For the canonical system of generators of the
algebra $A$, set:
\begin{equation} \label{eq:hom}
\varphi (b_\beta  ) = b_\beta   \cdot (1 + f), \  \varphi
(a_\alpha ) = a_a  + h_\alpha,  \quad \alpha  \in \Lambda,
\end{equation}
where here $h_\alpha \in A^2 $, $h_\alpha   = a_\alpha  a_{j_1 }
f_1 + \cdots + a_\alpha a_{j_q } f_q $. This defines the images of
the canonical set of generators of  $A$  under $\varphi$. Extend
the action of  $\varphi $ to $A$ agreeing with the definition of a
homomorphism.

\begin{prop} \label{prop:224}
 In this notation, the map $\varphi :A
\rightarrow A$ is a homomorphism. Moreover if the ideal $I$ is a
torsion-free module over $R$ then $\varphi $ is injective.
\end{prop}
\proof We check that Equation (\ref{eq:hom}) gives rise to a
homomorphism on $A$. We use the presentation of $A$ given in
Theorem \ref{thm:223}. Every element $b \in I$ can be written in
the form
$$
b = b_{\beta _1 } f^{(1)}  + \cdots + b_{\beta _l } f^{(l)},
$$
where $f^{(k)}  \in R$. Using induction on the degree of the
polynomials $f^{(k)}$ one verifies that $\varphi (b) = b \cdot (1
+ f)$:
$$
\varphi (b_\beta  a_\alpha  ) = \varphi (b_\beta  )\varphi
(a_\alpha  ) = (b_\beta  (1 + f))(a_\alpha   + h_\alpha  ) =
b_\beta  a_\alpha  (1 + f).
$$
Now it is fairly obvious that $\varphi $ preserves the relations
of $A$ of the first type (the ones that reflect the $R$-module
structure of $I$) and of the second type (the ones that express
the fact that $I$ is an Abelian ideal). The relations of the third
type take the form:
$$
a_i a_j  = b_{\beta _1 } f^{(1)} _{ij}  + \cdots + b_{\beta _l }
f^{(l)} _{ij},
$$
where $f_{ij} ^{(k)}  \in R$, $i,j \in \Lambda $ is an arbitrary
pair of indices. To verify that $\varphi $ preserves these
relations it suffices to check that $\varphi (a_i )\varphi (a_j )
= a_i a_j (1 + f)$:
\begin{gather} \notag
\begin{split}
\varphi (a_i )\varphi (a_j )  = (a_i + a_i a_{j_1 } f_1 & + \cdots
+ a_i a_{j_q } f_q )(a_j  + a_j a_{j_1 } f_1  + \cdots + a_j
a_{j_q} f_q ) = \\
 = a_i a_j  + (a_i a_{j_1 } a_j f_1  &+ \cdots + a_i
a_{j_q} a_j f_q ) - (a_j a_{j_1 } a_i f_1  + \cdots + a_j a_{j_q} a_i f_q ) = \\
= a_i a_j  + a_i a_j a_{j_1 } f_1  + \cdots +& a_i a_j a_{j_q} f_q
= a_i a_j (1 + f).
\end{split}
\end{gather}

Since the action of $\varphi$ on the ideal $I$ is the
multiplication of elements of the module by the polynomial $(1 +
f)$. The homomorphism $\varphi$ (in the event that $I$ is a
torsion-free $R$-module), is injective. Any element $a \in A$ can
be written in the form $a = c + b$, where $c \in V$ , $b \in I$.
If $c \ne 0$ then, obviously, $\varphi (a) \ne 0$ and if $c = 0$,
then $a = b$, i.e. $a \in I$. Consequently, in any case the
homomorphism $\varphi $ is injective on $A$. \hfill $\blacksquare$

\subsection{Free Metabelian Lie Algebras}

In this section we list some of the properties of free metabelian
Lie algebras and find defining relations of its Fitting's radical
treated as a module over the ring of polynomials $R$.

Let $F$ be the free metabelian Lie algebra over the field $k$. And
let $\left\{ a_\alpha  ,\,\alpha  \in \Lambda \right\}$ be a free
base of $\F$. Suppose that the set of indices $\Lambda $ is
totally ordered. Then left-normed monomials $a_{i_1 } a_{i_2 }
\cdots a_{i_m } $ satisfying $i_1  > i_2  \le i_3 \le \ldots \le
i_m $ are termed normalised. In \cite{Art1}  V. A. Artamonov
proves the following

\begin{thm} \label{thm:231}
 The set of all normalised monomials of the free metabelian Lie algebra
$\F$ forms its linear basis over $k$.
\end{thm}

\begin{cor} \label{cor:II-propofsubF} \quad \quad \quad
\begin{itemize}
    \item The images of  the free base $\left\{ a_\alpha,\,\alpha \in \Lambda \right\}$
            of the free Lie algebra $\F$ in the factor-algebra
            ${\raise0.7ex\hbox{$\F$} \!\mathord{\left/
            {\vphantom {\F {\F^2 }}}\right.\kern-\nulldelimiterspace}
            \!\lower0.7ex\hbox{${\F^2 }$}}$ form its additive basis over $k$.
    \item  If $|\Lambda |> 1$ then the algebra  $\F$ is not
    Abelian. And in the event that $|\Lambda | = 1$ the algebra  $\F$
is a one-dimensional vector space.
    \item Every collection of linearly independent modulo $\F^2 $ elements $\left\{ b_\beta
,\,\beta \in \Beta \right\} $ of the algebra $\F$ generates a free
metabelian Lie algebra of the rank $|\Beta|$.
\end{itemize}
\end{cor}

The commutant of an arbitrary  metabelian Lie algebra lies in its
Fitting's radical. It turns out to be that if $F$ is non-Abelian
then this inclusion is in fact an equality.

\begin{prop}\label{prop:fit=f2} \label{prop:22}
If  $F$ is a non-Abelian free metabelian Lie algebra then \/ {\rm
$\Fit(F) = F^2$}.
\end{prop}
\proof Take an arbitrary element $a \notin F^2 $ and show that $a
\notin \Fit(F)$. Write $a$ as a linear combination of normalised
words:
$$
a = \alpha _1 a_{i_1 }  + \cdots + \alpha _n a_{i_n }  + b,
$$
where $b \in F^2 $. Since $a \notin F^2 $, we have $c = \alpha _1
a_{i_1 }  + \cdots + \alpha _n a_{i_n }  \notin F^2 $. It suffices
to show that $c \notin \Fit(F)$. Assume the converse: $c \in
\Fit(F)$. The element $c$, therefore, can be embodied into a free
base of the algebra $F$. The new base is obtained from the initial
one by a linear transformation. Since $c \in \Fit(F)$, the element
$c$ lies in a nilpotent  ideal $I$ of $F$. For $F$ is non-Abelian,
there exists an element $d$ from the free base of $F$ so that $d
\circ c \ne 0$, $dc \in I$, where $c \in I$. Consequently,
multiple multiplication of $dc$ by $c$ results zero, i.e.
$dccc\cdots c = 0$. This derives a contradiction with Theorem
\ref{thm:231}. \hfill $\blacksquare$

On account of Proposition \ref{prop:fit=f2}, the free base
$\left\{ a_\alpha ,\,\alpha  \in \Lambda \right\}$ of a free
non-Abelian Lie algebra $F$ is a maximal linearly independent
modulo $\Fit(F)$ system of elements. The system of generators of
$\Fit(F)$, treated as a module over $R = k[x_\alpha  ,\alpha  \in
\Lambda ]$,
 (denote it $C(F)$) is constructed as follows.
The set  $\left\{ a_\alpha  a_\beta \mid {\;\alpha ,\beta  \in
\Lambda } \right\}$ generates the $R$-module $\Fit(F)$. For each
pair of indices $(\alpha ,\beta )$ either $a_\alpha a_\beta  $ or
$a_\beta  a_\alpha$ only, is included in $C(F)$. Consequently, the
canonical system of generators of the algebra $F$ is the union of
$\left\{ a_\alpha  ,\,\alpha  \in \Lambda \right\}$ and $C(F)$.
According to the Jacoby identity, the generators of  $\Fit(F)$
satisfy relations of the type:
\begin{equation} \label{eq:relforfit}
a_\alpha  a_\beta  x_\gamma   + a_\beta  a_\gamma  x_\alpha   +
a_\gamma  a_\alpha  x_\beta   = 0.
\end{equation}
Lemma \ref{lem:relforfit} below shows that Relations
(\ref{eq:relforfit}) define the module $\Fit(F)$.

\begin{lem} \label{lem:relforfit}
The system of generators $C(F)$  and the set of Relations {\rm
(\ref{eq:relforfit})} given above define an $R$-module
presentation of the Fitting's radical \/ {\rm $\Fit(F)$} of the
free metabelian Lie algebra $F$.
\end{lem}
\proof Let $M$ be the $R$-module given by the generators and
relations described above. Consider the following formal direct
sum $V \oplus M$, where $V$ is a $k$-vector space with the base
$\left\{ a_\alpha  ,\,\alpha  \in \Lambda \right\}$. The set $V
\oplus M$ is a $k$-vector space in the obvious way. We define a
multiplication in $V \oplus M$ by setting:
$$
(a_\alpha   \oplus m_1 ) \circ (a_\beta   \oplus m_2 ) = (0 \oplus
a_\alpha  a_\beta   + m_1 x_\beta   - m_2 x_\alpha  ).
$$
The reader can check that the space $V \oplus M$ with the above
operations is a metabelian Lie algebra with the generators
$\left\{ a_\alpha ,\,\alpha \in \Lambda \right\}$. The commutant
of $V \oplus M$ is the module  $M$. Since  $F$ is a free algebra,
there exists the canonical epimorphism  $\varphi :F \rightarrow V
\oplus M$. The image of the commutant of $F$ coincides with $M$
and the statement of the lemma follows. \hfill $\blacksquare$

\section{$U$-algebras} \label{sec:3}
In this section we construct special matrix metabelian Lie
algebras.

Every matrix metabelian Lie algebra is constructed by a ring of
polynomials and a free module $T$ over $R$. In Section \ref{ss:31}
we give a description of these algebras and prove that free
metabelian Lie algebra of an arbitrary rank can be embedded into a
matrix metabelian Lie algebra. The idea of this embedding arises
from papers  \cite{Art1} and \cite{Art2}.

 In Section \ref{ss:32}
we introduce the notion of a $U$-algebra and establish connections
between $U$-algebras and matrix metabelian Lie algebras. In
particular, every finitely generated subalgebra of a $U$-algebra
also embeds into a special matrix algebra.

Matrix metabelian Lie algebras, therefore, provide a convenient
presentation of free metabelian Lie algebras and $U$-algebras. We
shall use this presentation further, when instead of working with
metabelian Lie algebras we investigate commutative rings and
modules.

\subsection{Matrix Metabelian Lie Algebras} \label{ss:31}

Let $k$ be a field, let $R$ be the ring of polynomials from
commuting variables $\left\{ x_\alpha  ,\alpha  \in \Lambda
\right\}$ with coefficients from  $k$ and let $T$ be the free
module over the ring $R$ with the base  $\left\{ u_i ,i \in
I\right\}$. By this ring $R$ and module $T$ we construct a matrix
metabelian Lie algebra. Denote by $M_{I,\Lambda}$ the set of all
matrices of the form
     $\left( {\begin{array}{*{20}c}
        f & u  \\
        0 & 0  \\
     \end{array}} \right)$,
where $f \in R$ and $u \in T$. We turn the set $M_{I,\Lambda } $
into a Lie $k$-algebra by setting:
\begin{gather} \notag
\begin{split}
\alpha  \cdot \left( {\begin{array}{*{20}c}
        f & u  \\
        0 & 0  \\
     \end{array}} \right) = \left( {\begin{array}{*{20}c}
        {\alpha  \cdot f} & {\alpha  \cdot u}  \\
        0 & 0  \\
     \end{array}} \right), \alpha  \in k; \\
      \left( {\begin{array}{*{20}c}
        f & u  \\
        0 & 0  \\
     \end{array}} \right) + \left( {\begin{array}{*{20}c}
        g & v  \\
        0 & 0  \\
     \end{array}} \right) = & \left( {\begin{array}{*{20}c}
        {f + g} & {u + v}  \\
        0 & 0  \\
     \end{array}} \right); \\
     \left( {\begin{array}{*{20}c}
        f & u  \\
        0 & 0  \\
     \end{array}} \right) \circ & \left( {\begin{array}{*{20}c}
        g & v  \\
        0 & 0  \\
     \end{array}} \right) = \left( {\begin{array}{*{20}c}
        0 & {ug - vf}  \\
        0 & 0  \\
     \end{array}} \right).
\end{split}
\end{gather}

The set $M_{I,\Lambda}$ is closed under the operations introduced
above and forms a metbabelian Lie algebra. The correctness of the
anti-commutativity identity, the metabelian identity, the axioms
of a vector space and the distributivity of multiplication in
$M_{I,\Lambda } $ is straightforward. We next show that
$M_{I,\Lambda}$ satisfies the Jacoby identity. For the sake of
brevity we write elements from $M_{I,\Lambda}$ as pairs: $\left\{
(f,u)\right\}$.

The product of arbitrary three elements from  $M_{I,\Lambda } $
equals:
$$
(f_1 ,v_1 )\circ (f_2 ,v_2 ) \circ (f_3 ,v_3 ) = (0,v_1 f_2  - v_2
f_1 ) \circ (f_3 ,v_3 ) = (0,v_1 f_2 f_3  - v_2 f_1 f_3 ).
$$

So the Jacoby identity takes the form:
\newline $(f_1 ,v_1 ) \circ (f_2 ,v_2 ) \circ (f_3
,v_3 ) + (f_2 ,v_2 ) \circ (f_3 ,v_3 ) \circ (f_1 ,v_1 ) + (f_3
,v_3 ) \circ (f_1 ,v_1 ) \circ (f_2 ,v_2 )  = \\= (0,v_1 f_2 f_3 -
v_2 f_1 f_3 ) + (0,v_2 f_3 f_1  - v_3 f_2 f_1 ) + (0,v_3 f_1 f_2 -
v_1 f_3 f_2 ) = (0,0).$

We thereby have proven that the set of matrices $M_{I,\Lambda } $
with the operations above  is a metabelian Lie $k$-algebra. Such
algebras and all their subalgebras are called \emph{matrix
metabelian Lie algebras}.

We allow for the possibility that $\Lambda = \emptyset$, i.e. $R =
k$, in which case $T$ is a vector space over $k$.

In general case, the algebra $M_{I,\Lambda}$ is not Abelian:
$$
(0,u) \circ (f,u) = (0,uf) \ne (0,0),
$$
provided that $f \ne 0$ and $u \ne 0$. Although $M_{I,\Lambda } $
contains Abelian subalgebras. For instance, the commutant or any
one-generated subalgebra are Abelian. The multiplication in
$M_{I,\Lambda } $ is defined so that all the elements (matrices)
of its commutant have zero principal diagonal,  i.e. the commutant
consists of elements of the form $(0,u)$. It turns out to be that
the collection of all such elements forms the Fitting's radical
$\Fit(M_{I,\Lambda } )$. The foregoing discussion is summarised
and generalised in the lemma below.

\begin{lem} \label{lem:311}
Let $A$ be a non-Abelian subalgebra of matrix metabelian Lie
algebra $M_{I,\Lambda } $.  Then the  Fitting's radical of $A$
consists of all the elements from  $A$ with zero principal
diagonal:{\rm
$$
        \Fit(A) = \left\{ {(f,u) \in A,\;\,f = 0} \right\}.
$$}
\end{lem}
\proof Let $B$ be the following set $\left\{ {(f,u) \in A,\;\;f =
0} \right\}$. We show that  $B = \Fit(A)$. Take an arbitrary
element $(0,u) \in B$ and let $I$ be the principal ideal of $A$
generated by $(0,u)$. Clearly $I \subset B$, and so $I$ is an
Abelian ideal and, thus  $(0,u) \in \Fit(A)$. To show the reverse
inclusion take  an arbitrary element $z = (f,u)$ from $A$. Suppose
that $z\notin B$, in particular $f \ne 0$. Assume the converse: $z
\in \Fit(A)$. Since $z \in \Fit(A)$, by Corollary
\ref{cor:II-fitchar}, $z$ lies in a nilpotent ideal $I$. Since $A$
is non-Abelian, there exist two elements in $A$ such that their
product is non-zero. By the definition, the product of two
arbitrary elements lies in $B$. Consequently, there is to exist a
non-zero element $w = (0,v)\in B$, $v \ne 0$. Then, in the above
notation, $w \circ z = (0,vf) \in I$ and $w \circ z \ne 0$. For
$I$ is a nilpotent ideal, there exists a positive integer $n$ such
that  $w \circ z \circ \cdots \circ z = (0,vf^n ) = (0,0)$, what
derives a contradiction. \hfill $\blacksquare$

\begin{rem}
There exists an Abelian subalgebra  $A$ of the matrix metabelian
Lie algebra $M_{I,\Lambda}$ so that Lemma {\rm \ref{lem:311}}
fails. For instance, let $A$ be the algebra generated by  $(f,u)$,
$f \ne 0$. In that case {\rm $\Fit(A) = A$}.
\end{rem}

\begin{cor}\label{cor:31}
The Fitting's radical of an arbitrary subalgebra of a matrix
metabelian Lie algebra is Abelian.
\end{cor}
\textbf{Notation} \textit{Let $M_{I,\Lambda } $ be a matrix
metabelian Lie algebra. By $M_{I,\Lambda }^0 $ denote the
following subset of $M_{I,\Lambda }$:
 $$
 M_{I,\Lambda }^0  = \{ \left. {(f,u)} \right|\;f = \beta _{\alpha _1 } x_{\alpha _1 }  + \beta _{\alpha _2 } x_{\alpha _2 }
 + \cdots + \beta _{\alpha _m } x_{\alpha _m } ,\;\beta _{\alpha _j }  \in k\} .
 $$}

Clearly, $M_{I,\Lambda }^0 $ is a subalgebra of  $M_{I,\Lambda }$.
In the following we shall be particularly interested in a more
narrow class of matrix metabelian Lie algebras, namely, the
subalgebras of $M_{I,\Lambda }^0 $.

\begin{lem} \label{lem:312}
Let $A$ be a subalgebra of a matrix metabelian Lie algebra
$M_{I,\Lambda }^0 $. Then the Fitting's radical {\rm $\Fit(A)$} is
a torsion-free module over the ring of polynomials.
\end{lem}
\proof On behalf of  Corollary \ref{cor:31}, $\Fit(A)$ is Abelian.
Consequently, it allows a structure of a module over the ring of
polynomials. Let $\left\{ e_d ,\,d \in \Delta \right\}$ be a basis
of the vector space ${\raise0.7ex\hbox{$A$} \!\mathord{\left/
{\vphantom {A {\Fit(A)}}}\right.\kern-\nulldelimiterspace}
 \!\lower0.7ex\hbox{${\Fit(A)}$}}$, let $\left\{ {(f_d ,u_d ),\;d \in \Delta } \right\}$
be the set of preimages in  $A$ of elements of this basis. All the
polynomials $f_d $'s are $k$-linear combinations of letters
$\left\{ x_\alpha  ,\;\alpha  \in \Lambda \right\} $. Furthermore,
on account of Lemma \ref{lem:311}, the set of polynomials $\left\{
{f_d ,\;d \in \Delta } \right\}$ is linearly independent.

We need to show that $\Fit(A)$ is a torsion-free module over the
ring of polynomials $R = k\left[x_d ,\,d \in \Delta \right]$. Take
an arbitrary element $b \in \Fit(A)$, a non-zero polynomial $g(x_d
) \in R$ and consider the product $b \cdot g$. By  Lemma
\ref{lem:311} every element of $\Fit(A)$ has the form $b = (0,v)$
and, consequently, $b \cdot g = (0,vg(f_d ))$. The polynomial
$g(f_d )$ is also non-zero and therefore $b \cdot g \ne 0$. \hfill
$\blacksquare$

We  next show how a free metabelian Lie algebra $F$  over the
field $k$ is embedded into a matrix algebra. Let $\left\{ a_\alpha
,\,\alpha  \in \Lambda \right\}$ be the free base of $F$. Then the
algebra $F$ is a subalgebra of the matrix metabelian Lie algebra
$M_{\Lambda ,\Lambda }^0 $. The matrix algebra $M_{\Lambda
,\Lambda }^0 $ is constructed by the ring of polynomials $R =
k\left[x_\alpha  ,\alpha  \in \Lambda \right]$ and by the free
module $T$ with the set of free generators $\left\{ u_\alpha
,\alpha  \in \Lambda \right\}$. Consider a subalgebra $L$ in
$M_{\Lambda ,\Lambda }^0 $ generated by
$$
    \left\{ {z_\alpha   = \left( {\begin{array}{*{20}c}
        {x_\alpha  } & {u_\alpha  }  \\
        0 & 0  \\
     \end{array}} \right),\;\alpha  \in \Lambda } \right\}.
$$
The cardinality of the set $\left\{ {z_\alpha ,\;\alpha \in
\Lambda } \right\}$ of generators  of the algebra $L$ coincides
with the cardinality of the free base $\{ a_\alpha ,\,\alpha \in
\Lambda \} $ of $F$. Consequently, there exists the canonical
homomorphism $\gamma :F \rightarrow M_{\Lambda ,\Lambda }^0 $ so
that $\gamma (a_\alpha  ) = z_\alpha ,\alpha \in \Lambda $. The
following theorem shows that the algebra $L$ is a presentation of
$F$ in a matrix algebra.

\begin{thm} \label{thm:hominmtx}
In the above notation, the homomorphism $\gamma $ is an embedding
of the free metabelian Lie algebra $F$ into the matrix metabelian
Lie algebra $M_{\Lambda ,\Lambda }^0 $.
\end{thm}
\proof Let $I$ be the kernel of $\gamma $. We prove that $I$ is
the trivial ideal.

By Theorem \ref{thm:231} the algebra $L$ (as a $k$-vector space)
is generated by normalised words from the letters $\left\{
{z_\alpha ,\;\alpha  \in \Lambda } \right\}$. Recall that
normalised words are left-normed words $z_{i_1 } z_{i_2 } \cdots
z_{i_m } $ such that $i_1  > i_2  \le i_3  \le \ldots \le i_m$ (we
assume that the set of indices $\Lambda $ is totally ordered). To
prove that $I$ is the trivial ideal it suffices to show that the
normalised words in $L$ are linearly independent.

In the matrix representation the normalised word $z_{i_1 } z_{i_2
}\cdots z_{i_m } $ is presented by the following matrix
$$
    \left( {\begin{array}{*{20}c}
        0 & {(u_{i_1 } x_{i_2 }  - u_{i_2 } x_{i_1 } ) \cdot x_{i_m }  \cdots x_{i_3 } }  \\
        0 & 0  \\
     \end{array}} \right).
$$
Since $(u_{i_1 } x_{i_2 }  - u_{i_2 } x_{i_1 } ) \cdot x_{i_m }
\cdots x_{i_3 } $ is a non-zero element of the free module $T$,
this matrix is non-zero. Assume further that
\begin{equation} \label{eq:n1}
\gamma (\sum\limits_{\bar i} {\alpha _{\bar i} a_{i_1 } a_{i_2 }
\cdots a_{i_m } } ) = 0,\;\;\alpha _{\bar i}  \in k,
\end{equation}
and show that for all multi-indices $\bar i$ the coefficient
$\alpha _{\bar i}  = 0$. Obviously, we may assume that for all
$\bar i$ in Equation (\ref{eq:n1})  $m \ge 2$. Equation
(\ref{eq:n1}) implies the following equality in the module
 $T$:
\begin{equation} \label{eq:n2}
    \sum\limits_{\bar i} {\alpha _{\bar i} (u_{i_1 } x_{i_2 }  -
    u_{i_2 } x_{i_1 } ) \cdot x_{i_k }  \cdots x_{i_3 } }  =
    0.
\end{equation}
Take the maximal index $i_1$  among all first indices in all
multi-indices $\bar i$ involved in  Equation (\ref{eq:n1}). By the
definition of a normalised word, we obtain that the coefficient of
$u_{i_1 }$ on the left-hand side of Equation (\ref{eq:n2}) is
$$
\sum\limits_{\bar i = i_1 ,\ldots,i_m } {\alpha _{\bar i} x_{i_m }
\cdots x_{i_3 }  \cdot x_{i_2 } }.
$$
Since $T$ is a free module and since $u_{i_1 } $ is an element
from the basis of  $T$, this coefficient equals zero in the ring
of polynomials $R$. Again, by the definition of a normalised word,
we obtain that if all the $\alpha _{\bar i}$'s are zero, provided
that the multi-index $\bar i$ begins with $i_1 $. We proceed by
induction on the first coordinate of the multi-index $\bar i$ and
the statement follows. \hfill $\blacksquare$

As a conclusion we formulate and prove several simple and
important properties of elements of matrix metabelian Lie
algebras.

\begin{lem} \label{lem:forax}
Let $M_{I,\Lambda } $ be an arbitrary matrix metabelian Lie
algebra and let $x = (f_1 ,u_1 )$, $y = (f_2 ,u_2 )$, $z = (f_3
,u_3 )$ be elements from $M_{I,\Lambda}$. Then
    \begin{enumerate}
        \item If $xyx = 0$ and $xyy = 0$ then $xy = 0$.
        \item If $xy = 0$, $xz = 0$ and $x \ne 0$, then $yz = 0$.
    \end{enumerate}
\end{lem}
\proof Let $x,y \in M_{I,\Lambda }$, $xyx = 0$ and $xyy = 0$. In
terms of matrices two latter conditions rewrite as follows:
$$(0,(u_1 f_2  - u_2 f_1 )f_1 ) = (0,0) \hbox{ and } (0,(u_1 f_2  - u_2
f_1 )f_2 ) = (0,0).
$$
Now the required equality is essentially
immediate: $xy = (0,u_1 f_2 - u_2 f_1 ) = (0,0)$.

We now prove the second statement. The conditions imposed imply:
\begin{equation} \label{eq:m1}
u_1 f_2  = u_2 f_1,
\end{equation}
\begin{equation} \label{eq:m2}
        u_1 f_3  = u_3 f_1.
\end{equation}
Furthermore, $f_1  \ne 0$ and $u_1  \ne 0$. We need to show that
$yz = 0$, i.e.
\begin{equation} \label{eq:m3}
        u_2 f_3  = u_3 f_2.
\end{equation}
In the event that $f_1  \ne 0$ Equation (\ref{eq:m3}) is
equivalent to the following equality $u_2 f_3 f_1 = u_3 f_2 f_1$.
Multiplying Equations (\ref{eq:m1}) and (\ref{eq:m2}) by $f_3 $
and $f_2 $ correspondingly we obtain $u_2 f_3 f_1 = u_3 f_2 f_1$.
Now, if $f_1  = 0$ then $f_2  = f_3  = 0$ and Equation
(\ref{eq:m3}) also holds. \hfill $\blacksquare$

\begin{cor}
According to Theorem {\rm \ref{thm:hominmtx}}, Lemma {\rm
\ref{lem:forax}} holds in the free metabelian Lie algebra.
\end{cor}

\subsection{$U$-Algebras} \label{ss:32}

\begin{defn}
We term a metabelian Lie algebra $A$ over a field $k$ a metabelian
Lie $U$-algebra if
\begin{itemize}
    \item {\rm $\Fit(A)$} is Abelian;
    \item {\rm $\Fit(A)$}, treated as a module over the ring of
    polynomials (defined as in Section \ref{sec:21}), is torsion-free.
\end{itemize}
\end{defn}

For the sake of brevity, we sometimes call metabelian Lie $U$-algebras by $U$-algebras.

\begin{rem} The Fitting's radical \/ {\rm $\Fit(A)$} admits the structure of
a module over the ring  $R = k\left[ x_\alpha  ,\alpha  \in
\Lambda \right]$. The system of linearly independent modulo {\rm
$\Fit(A)$} elements $\left\{ a_\alpha ,\alpha \in \Lambda \right\}
$ of $A$ may vary depending on the definition of a module
structure on {\rm $\Fit(A)$}. The transformation between two such
systems is a $k$-linear change of variables of the ring $R$. This
implies that if for $\left\{ a_\alpha ,\alpha \in \Lambda
\right\}$ the module {\rm $\Fit(A)$} is  torsion-free then it is
torsion-free for another such tuple of elements. Therefore the
definition of a $U$-algebra is correct.
\end{rem}

If $A$ is an Abelian Lie algebra then $A = \Fit(A)$, thus the
Fitting's radical $\Fit(A)$ is a module over the field $k$.
Consequently, every Abelian Lie algebra is a  $U$-algebra.

The connection between metabelian $U$-algebras and matrix
metabelian Lie algebras is established by the following

\begin{thm} \label{thm:321}
Every finitely generated metabelian Lie $U$-algebra is isomorphic
to a subalgebra of a matrix metabelian Lie algebra $M_{I,\Lambda
}^0$. Conversely,  every subalgebra of the algebra $M_{I,\Lambda
}^0 $ is a $U$-algebra.
\end{thm}
\proof Directly from Lemmas \ref{lem:311} and \ref{lem:312} we
obtain that every subalgebra $A$ of the matrix metabelian Lie
algebra $M_{I,\Lambda }^0 $ is a $U$-algebra.

Let $A$ be a finitely generated $U$-algebra. By the definition
$\Fit(A)$ is Abelian and, therefore, admits the structure of a
module over the ring of polynomials. Since $A$ is a finitely
generated algebra, by Lemma \ref{lem:222}, its canonical system of
generators $a_1 ,\ldots,a_r$, $b_1 ,\ldots,b_l$ is finite. Again,
by the definition, $\Fit(A)$ is a torsion-free module over the
ring $R = k[x_1 ,\ldots, x_r ]$.

Any torsion-free finitely generated module over the ring of
polynomials embeds into a free module, thus $\Fit(A)$ embeds into
the free module  $T$ over the ring $R$ of the rank $s$, $s \le l$
(see \cite{Bourb, Leng}). Let $\varphi :\Fit(A) \rightarrow T$ be
the correspondent embedding. Using this embedding we construct an
embedding of the algebra $A$, $\gamma :A \rightarrow M_{s,r}^0 $.

By the definition of $\gamma$, for the standard generators we set:
\begin{gather} \notag
\begin{split}
\gamma (b_i )& = (0,\varphi (b_i ) \cdot (\sum\limits_{k = 1}^r
{x_k \,} )),\;i = 1,\ldots,l, \\
&\gamma (a_j ) = (x_j ,\sum\limits_{k = 1}^r {\varphi (a_j  \circ
a_k )} ),\;j = 1,\ldots,r.
\end{split}
\end{gather}
 We extend the definition of $\gamma $ on $A$ following the
definition of a homomorphism. We are left to check that $\gamma $
is a correct homomorphism.

To prove that $\gamma $ is a correct homomorphism we prove that it
preserves the three types of defining relations of $A$ given in
Theorem \ref{thm:223}.

The first type of defining relations of $A$ are exactly the
defining relations of $\Fit(A)$ treated as a module over the ring
$R$. Every module relation of $\Fit(A)$ has the form: $ b_1  \cdot
f_1 + \cdots + b_l  \cdot f_l  = 0$, $f_i  \in R$. For every pair
of indices $(i,j)$, $i = 1,\ldots, l$ and
 $j = 1,\ldots,r$ we obtain:
$$
\gamma (b_i  \circ a_j ) = \gamma (b_i ) \circ \gamma (a_j ) = (0,\varphi (b_i ) \cdot (\sum\limits_{k = 1}^r {x_k \,} ) \cdot x_j ).
$$
Now for an arbitrary polynomial $f \in R$, applying induction on
the degree of a polynomial, we obtain:
$$
\gamma (b_i  \cdot f) = (0,\varphi (b_i ) \cdot (\sum\limits_{k = 1}^r {x_k \,} ) \cdot f) = (0,\varphi (b_i  \cdot f) \cdot (\sum\limits_{k = 1}^r {x_k \,} )).
$$
This gives the required equality:
\begin{gather} \notag
\begin{split}
\gamma (b_1  \cdot f_1  + \cdots + b_l  \cdot f_l ) = \gamma (b_1
\cdot f_1 ) + \cdots + \gamma (b_l  \cdot f_l ) = \\
 = (0,\varphi (b_1  \cdot f_1  + \cdots + b_l  \cdot f_l ) \cdot &(\sum\limits_{k = 1}^r {x_k } )) = (0,0).
\end{split}
\end{gather}
By Lemma \ref{lem:311}, $\gamma (b_i ) \in \Fit(M_{s,r}^0 )$, $i =
1,\ldots,l$.

Consequently, the relations of the second type, which express that
 $\Fit(A)$ is Abelian, are preserved under the action of $\gamma $.

The third type of relations of the algebra $A$ shows that $a_i
\circ a_j  \in \Fit(A)$ for every pair of indices $(i,j)$. These
relations have the form: $a_i \circ a_j = b_1 \cdot g_{ij} ^1  +
\cdots + b_l  \cdot g_{ij} ^l ,\;\,g_{ij} ^p \in R$. We have
already shown above that
$$
\gamma (b_1  \cdot g_{ij} ^1  + \cdots + b_l  \cdot g_{ij} ^l ) =
(0,\varphi (b_1  \cdot g_{ij} ^1  + \cdots + b_l  \cdot g_{ij} ^l
) \cdot (\sum\limits_{k = 1}^r {x_k } )).
$$
We therefore are left to show that $\gamma (a_i \circ a_j ) =
(0,\varphi (a_i  \circ a_j ) \cdot (\sum\limits_{k = 1}^r {x_k }
))$. By the means of direct computation
\begin{gather} \notag
\begin{split}
\gamma (a_i  \circ a_j ) = \gamma (a_i ) \circ \gamma (a_j ) =
(x_i ,\sum\limits_{k = 1}^r {\varphi (a_i  \circ a_k )} ) \circ
(x_j ,\sum\limits_{k = 1}^r {\varphi (a_j  \circ a_k )} )
=\\
= (0,\sum\limits_{k = 1}^r {(\varphi (a_i  \circ a_k ) \cdot x_j -
\varphi (a_j  \circ a_k ) \cdot x_i )} ) = (0,\varphi (a_i \circ
a_j ) \cdot (\sum\limits_{k = 1}^r &{x_k } )).
\end{split}
\end{gather}
The latter equality follows from the Jacoby identity.

To prove the theorem it suffices to show that the kernel of the
homomorphism $\gamma $ is trivial. Let $a \in A$. We show that
 $a \ne 0$ implies that $\gamma (a) \ne (0,0)$. The element $a$
 has a unique presentation in the form:
 $$
 a = \alpha _1 a_1  + \cdots + \alpha _n a_n  + b,
 $$
where $\alpha _i  \in k$, $b \in \Fit(A)$. If $\alpha _1 a_1 +
\cdots + \alpha _n a_n  \ne 0$ then $\gamma (a) \ne (0,0)$. Let $a
= b \in \Fit(A)$ and $b \ne 0$. The element $b$ can be written as
 $b = b_1 f^1  + \cdots + b_l f^l $, where $f^k  \in R$.  Now, since $\varphi $
is an embedding and $T$ is a torsion free module over $R$, we
conclude that
$$
\gamma (b) = (0,\varphi (b_1  \cdot f^1  + \cdots + b_l  \cdot f^l
) \cdot (\sum\limits_{k = 1}^r {x_k } )) \ne (0,0).
$$

We, therefore, have shown that $\gamma $ is an injective
homomorphism, whose restriction on the $U$-algebra $A$ is an
embedding of $A$ into the matrix metabelian Lie algebra
$M_{s,r}^0$ \hfill $\blacksquare$

\begin{rem} \label{rem:322}
In the extreme case when  $A$ is an Abelian algebra, $A$ embeds
into $M_{l,0}^0 $, where $l$ is the dimension of the $k$-vector
space $A$.
\end{rem}

\begin{cor}
Notice that free metabelian Lie algebras are $U$-algebras.
\end{cor}

Theorem \ref{thm:321} characterises finitely generated
$U$-algebras. Further results give a description of arbitrary
metabelian $U$-algebras.

\begin{lem} \label{lem:322}
Let  $A$ be a metabelian Lie $U$-algebra and let $B$ be a
subalgebra of  $A$. Then $B$ is a $U$-algebra.
\end{lem}
\proof If $B$ is an Abelian subalgebra of $A$, then, according to
Remark \ref{rem:322}, $B$ is a  $U$-algebra. We, therefore, assume
that $B$ is non-Abelian.

We show that $\Fit(B) = B \cap \Fit(A)$. The inclusion  $\Fit(B)
\supseteq B \cap \Fit(A)$ is obvious. The converse inclusion is at
issue. Take an arbitrary element $b \in B\smallsetminus \Fit(A)$.
Since $B$ is non-Abelian, there exist two elements $c,d \in B$ so
that $cd \ne 0$. Since $b \notin \Fit(A)$ and since $0 \ne cd \in
\Fit(A)$, all the monomials of the form $cdbb\cdots b$ are
non-zero. This implies that the element $b$ is contained in
neither of nilpotent ideals of the subalgebra  $B$ and
consequently $b \notin \Fit(B)$. This proves the reverse
inclusion.

It now follows that  $\Fit(B)$ is Abelian. Finally, any linearly
independent modulo $\Fit(B)$ system of elements from  $B$ is
linearly independent modulo $\Fit(A)$, thus $\Fit(B)$ is
torsion-free. \hfill $\blacksquare$

\begin{thm} \label{thm:323}
Let $A$ be a metabelian Lie  algebra. Then $A$ is a $U$-algebra if
and only if every finitely generated subalgebra of $A$ is a
$U$-algebra.
\end{thm}
\proof Let $A$ be a $U$-algebra. By Lemma \ref{lem:322} every
subalgebra of $A$ is a $U$-algebra. We are to prove the converse.
Suppose that every finitely generated subalgebra of $A$ is a
$U$-algebra.

We first show that $\Fit(A)$ is Abelian. Let $c,d \in \Fit(A)$ and
let $D =  \left< c,d \right> $ be the subalgebra generated by
$c,d$. Since $\Fit(D) \supseteq D \cap \Fit(A)$, both $c$ and $d$
lie in $\Fit(D)$. By the assumption  $D$ is a $U$-algebra and,
consequently, $cd = 0$.

We next show that $\Fit(A)$ is a torsion free module over the ring
of polynomials. Let $0 \ne b \in \Fit(A)$, and let $f(a_1
,\ldots,a_n )$ be a non-zero polynomial. Set $B =  \left< a_1
,\ldots,a_n ,b\right>$. The algebra $B$ is a finitely generated
subalgebra of  $A$ and thus is a $U$-algebra. Clearly,  $b \in
\Fit(B)$. This implies that to prove the required inequality $b
\cdot f \ne 0$ we need to show that the elements $a_1 ,\ldots,a_n
$ are linearly independent modulo $\Fit(B)$. Assume the converse.
Let $a = \alpha _1 a_1  + \cdots + \alpha _n a_n ,\;\alpha _i \in
k$ be a non-trivial linear combination and $a \in \Fit(B)$. We
prove that in this case $a \in \Fit(A)$, what derives a
contradiction with linear independence of the elements $a_1
,\ldots,a_n $ modulo $\Fit(A)$. Since $\Fit(B)$ is Abelian and $b
\in \Fit(B)$, we conclude that $ab = 0$. Take an arbitrary element
$c \in \Fit(A)$. Put  $C =  \left< a,b,c \right> $. The algebra
$C$ is a finitely generated subalgebra of $A$, thus by the
assumption $C$ is a $U$-algebra. Now, since $b,c \in \Fit(C)$, we
have that $b \ne 0$, $ab = 0$ and $bc = 0$. On account of Theorem
\ref{thm:321} the algebra $C$ embeds into a matrix metabelian Lie
algebra $M_{I,\Lambda }^0$. By Lemma \ref{lem:forax} we conclude
that $ac = 0$. Finally, on account of Lemma \ref{lem:214}, $a \in
\Fit(A)$. \hfill $\blacksquare$

\begin{cor} \label{cor:324}
Theorems {\rm \ref{thm:321}} and {\rm \ref{thm:323}} provide the
following chararcterisation of metabelian Lie $U$-algebras:
\begin{center}
A metabelian Lie algebra is a $U$-algebra if and only if every its
finitely generated subalgebra embeds into a special matrix
metabelian algebra of the form $M_{I,\Lambda }^0 $.
\end{center}
\end{cor}

\begin{thm}
Let $A$ be a $U$-algebra,  $x,y,z \in A $. Then
\begin{enumerate}
    \item if \/ $xyx = 0$ and \/ $xyy = 0$ then \/ $xy = 0$.
    \item if \/ $xy = 0$, $xz = 0$ and \/ $x \ne 0$ then \/ $yz = 0$.
\end{enumerate}
\end{thm}
\proof Let $C =  \left< x,y,z \right> $ be a finitely generated
subalgebra of the algebra $A$. Then, on behalf of Corollary
\ref{cor:324}, the algebra $C$ embeds into an algebra of the type
$M_{I,\Lambda }^0 $. Consequently, using Lemma \ref{lem:forax},
the statement follows.
 \hfill $\blacksquare$

\section{Universal Classes and Extensions of the Fitting's
Radical} \label{sec:4}

Let  $A$ be a $U$-algebra over a field  $k$. In this section we
describe two types of extensions of the Fitting's Radical of a
$U$-algebra $A$ over $k$. The first is the localisation of the
Fitting's radical of $A$ treated as an $R$-module by a maximal
ideal $\Delta \vartriangleleft R$. This algebra is called
$\Delta$-local. In Section \ref{ss:42} we investigate the
construction of $\Delta$-local algebras and establish its further
properties. The second type of extensions of the Fitting's radical
of $A$ is the direct module extension of the Fitting's Radical.

These constructions play an important role in our study of the
universal closure of the algebra $A$, which is crucial in
constructing algebraic geometry over Lie algebras (in fact, over
arbitrary algebraic systems). Section \ref{ss:41} holds
preliminary definitions and preliminary results on universal
closures, other universal classes and languages.

\subsection{Universal Classes} \label{ss:41}

The object of this section, which arises from papers \cite{alg1}
and \cite{alg2} is to introduce the counterparts to
group-theoretic notions from \cite{alg1} and \cite{alg2} in the
category of Lie $k$-algebras. We also formulate a number of
results for Lie algebras, which are complete counterparts to the
results of papers \cite{alg1} and \cite{alg2}. The proofs are
similar and therefore omitted.

The standard first order language $L$ of the theory of Lie
algebras over a fixed field $k$ consists of a symbol for
multiplication `$\circ$', a symbol for addition `$+$', a symbol
for substraction `$-$', a set of symbols $\left\{ k_\alpha| \;
\alpha \in k \right\}$ for multiplication of the elements of Lie
algebras on the coefficients from the field $k$ and a constant
symbol `$0$' which denotes zero.

For the purposes of algebraic geometry over a fixed Lie algebra
$A$ one is to study the category of $A$--Lie algebras, so it is
more convenient to use a bigger language $L_A$. The language of
the category of $A$--Lie algebras consists of the language $L$
together with the set of constant symbols enumerated by the
elements of $A$
$$
L_A  = L \cup \left\{ c_a |\; a \in A \right\}.
$$

 A Lie algebra $B$
over a field $k$ is called an $A$--Lie algebra if and only if it
contains  a designated copy of $A$, which we shall for most part
identify with $A$. It is clear that an $A$--Lie algebra $B$ can be
treated as a model of the language $L_A$. A homomorphism $\varphi$
from an $A$--Lie algebra $B_1$ to an $A$--Lie algebra $B_2$ is an
$A$-homomorphisms of Lie algebras if it is the identity on $A$.
The family of all $A$--Lie algebras together with the collection
of all $A$--homomorphisms form a category in the obvious way.

Set $\hom _A (B_1, B_2)$ to be the set of all $A$--homomorphisms
from the $A$-Lie algebra $B_1$ to the $A$-Lie algebra $B_2 $.

To an $A$--Lie algebra $B$ we link several model-theoretical
classes of $A$--Lie algebras. {\rm
\begin{itemize}
    \item The variety $A - \var(B)$ generated by
 $B$ is the class of all $A$--Lie algebras that satisfy all
the identities of the language $L_A$, satisfied by $B$.
    \item The quasivariety $A - \qvar(B)$ generated by  $B$ is the class of all $A$--Lie algebras
that satisfy all the quasi identities of the language $L_A$,
satisfied by $B$.
    \item The universal closure $A - \ucl(B)$ generated by
$B$ is the class of all $A$--Lie algebras that satisfy all the
universal sentences of the language $L_A$, satisfied by $B$.
\end{itemize}
}

Here are the definitions of the identity, the quasi identity and
the universal sentence in the language  $L_A $. {\rm
\begin{itemize}
    \item An $A$--universal sentence of the language $L_A$  is a
            formula of the type
            $$
            \forall x_1 \ldots\forall x_n \,(\mathop  \bigvee \limits_{j =
            1}^s \mathop  \bigwedge \limits_{i = 1}^t (u_{ij} (\bar x,\bar
            a_{ij} ) = 0\; \wedge \;w_{ij} (\bar x,\bar c_{ij} ) \ne 0)),
            $$
            where $\bar x = (x_1 ,\ldots,x_n )$ is an $n$-tuple of variables,
            $\bar a_{ij} $ and $\bar c_{ij} $ are the sets of constants from
            the algebra $A$ and $u_{ij}$, $w_{ij} $ are the terms in the
            language $L_A$ from the variables $x_1 ,\ldots,x_n$. In the event
            that an $A$--universal sentence involves no constants from the
            algebra $A$ this notion turns into standard notion of universal
            sentence in the language $L$.
    \item An $A$--identity of the language $L_A$ is the formula of the type
            $$
            \forall x_1 \ldots \forall x_n (\mathop  \bigwedge \limits_{i =
            1}^m \;r_i (\bar x, \bar a_{ij}) = 0),
            $$
            where $r_i (\bar x)$ are the terms in the language $L_A$ from the
            variables $x_1 ,\ldots,x_n$. In the event that an $A$--identity
            involves no constants from the algebra $A$ this notion turns into
            standard notion of identity of the language $L$.
    \item An $A$--quasi identity of the language $L_A$  is a formula of the type
            $$
            \forall x_1 \ldots\forall x_n (\mathop  \bigwedge \limits_{i =
            1}^m \;r_i (\bar x,\bar a_{ij}) = 0 \rightarrow s(\bar x, \bar b)
            = 0),
            $$
            where $r_i (\bar x)$ and $s(\bar x)$ are the terms. Coefficients
            free analogue is the notion of quasi identity.
\end{itemize}}

\begin{rem}
For our purposes we also need to consider the classes {\rm
$\var(B)$, $\qvar(B)$} and {\rm $\ucl(B)$}, which, by the
definition, are the variety, the quasivariety and the universal
closure generated by
 $B$ in the first order language $L$. The classes
{\rm $\var(B)$, $\qvar(B)$} and {\rm $\ucl(B)$} are also
particular cases of the classes {\rm $A - \var(B)$, $A -
\qvar(B)$} and {\rm $A - \ucl(B)$}, with $A = \left\{ 0\right\}$.
\end{rem}

\begin{rem}
For universal classes of $A$--Lie algebras the following sequence
of inclusions hold {\rm
$$
A - \ucl(B) \subseteq A - \qvar(B) \subseteq A - \var(B).
$$
}
The first inclusion is obvious, while the second one is implied by
the fact that every identity is equivalent to a conjunction of a
finite number of quasi identities. For instance, the identity
$\forall x_1 \ldots \forall x_n (\mathop \bigwedge \limits_{i =
1}^m \;r_i (\bar x) = 0)$ is equivalent to the following set of
$m$  quasi identities $\forall x_1 \ldots\forall x_n \forall y(y =
y \rightarrow r_i (\bar x) = 0)$.
\end{rem}

Let $B_1$ and $B_2$ be two $A$--Lie algebras. An $A$--Lie algebra
$B_2$ is termed $A$--discriminated by the $A$--Lie algebra $B_1$
if for every finite subset $\left\{ b_1 ,\ldots,b_m \right\}$ of
nonzero elements from the algebra $B_1$ there exists an
$A$--homomorphism $\varphi :B_1 \rightarrow B_2$ such that
$\varphi (b_i ) \ne 0$ for every $i = 1, \ldots, m$.

\begin{thm} \label{lem:I-dis-ucl} \label{thm:41}
Let $B_1$ and $B_2$ be two $A$--Lie algebras such that $B_1$ is
$A$--discriminated by the $A$--Lie algebra $B_2$. Then {\rm $B_1
\in A - \ucl(B_2)$}.
\end{thm}
\proof Recall that to prove that $B_1 \in A - \ucl(B_2)$ it
suffices to show that every finite submodel of the algebra $B_1$
$A$--embeds into the algebra $B_2$. This is obvious, since $B_1$
is $A$-discriminated by the algebra $B_2$. \hfill $\blacksquare$

\subsection{$\Delta$-Local Lie Algebras} \label{ss:II-deltalocalg}
\label{ss:42}

Let  $A$ be a $U$-algebra over a field  $k$. Let $\left\{
{z_\alpha ,\;\alpha  \in \Lambda } \right\}$ be a maximal set of
linearly independent modulo $\Fit (A)$ elements from $A$. Denote
by $V$  the linear span of $\left\{ {z_\alpha ,\;\alpha \in
\Lambda } \right\}$ over the field $k$.

Let $\Delta  =  \left< x_\alpha ,\alpha \in \Lambda  \right>$ be
the maximal ideal of the ring $R = k\left[x_\alpha ,\alpha  \in
\Lambda \right]$ generated by the variables $\left\{ x_\alpha
,\alpha \in \Lambda  \right\}$ .

Denote by  $R_\Delta  $ the localisation of the ring $R$ by
$\Delta $ and by $\Fit_\Delta  (A)$ the localisation of the module
$\Fit(A)$ by the ideal $\Delta $, i. e. the closure of $\Fit(A)$
under the action of the elements  of $R_\Delta$ (for definitions
see \cite{Bourb} and \cite{Leng}). Consider the direct sum $V
\oplus \Fit_\Delta (A)$ of vector spaces over $k$. We next define
a structure of an algebra on $V \oplus \Fit_\Delta  (A)$. By the
definition, the multiplication of the elements from $V$ is
inherited from $A$, the multiplication
 in $\Fit_\Delta  (A)$ is trivial. Set
$$
u \circ z_\alpha = u \cdot x_\alpha  ,\quad u \in \Fit_\Delta
(A),\;z_\alpha \in V,\quad u \cdot x_\alpha   \in \Fit_\Delta (A).
$$
and extend the definition of the multiplication by the elements
from $\Fit_\Delta (A)$ on the elements from $V$ to be linear.

One verifies that
\begin{itemize}
    \item this multiplication turns the vector space $V \oplus \Fit_\Delta
       (A)$ into a metabelian Lie algebra which we denote by
       $A_\Delta$,
    \item ${\raise0.7ex\hbox{${A_\Delta }$} \!\mathord{\left/
 {\vphantom {{A_\Delta  } {\Fit(A_\Delta  )}}}\right.\kern-\nulldelimiterspace}
\!\lower0.7ex\hbox{${\Fit(A_\Delta  )}$}} \cong V$,
    \item $\Fit(A_\Delta  )$ admits a structure of a module over
    $R$,
    \item $\Fit(A_\Delta  ) = \Fit_\Delta  (A)$,
    \item $A_\Delta  $ is a $U$-algebra,
    \item the algebra $A$ is a subalgebra of $A_\Delta$,
    \item even in the event that $A$ is finitely generated the algebra $A_\Delta  $
        is not finitely generated.
\end{itemize}

In the above notation, the algebra $A_\Delta$ is termed the
\emph{$\Delta$-localisation} of the algebra $A$.

We now point out some of the connections between the
$\Delta$-localisation and the universal closure of $A$.

\begin{lem} \label{lem:II-subofdel} \label{lem:421}
Let $A$ be a $U$-algebra and let $A_\Delta  $ be its
$\Delta$-localisation. Assume that $B$ is a finitely generated
subalgebra of $A_\Delta$. Then the algebra $B$ embeds into $A$.
\end{lem}
\proof Let $K = \left\{ b_1 ,\ldots,b_n \right\} $ be a set of
generators of a subalgebra of $B$. We next construct an injective
homomorphism  $\varphi: A_\Delta \rightarrow A_\Delta $ so that
$\varphi (b_i ) \in A$ ($A \subseteq A_\Delta  $) for each $b_i
\in K$. It is straightforward, that the restriction of  $\varphi $
to $B$ is an insertion of $B$ into $A$, $\varphi :B \rightarrow
A$. By the definition, the action of $\varphi $ is the
multiplication  by a polynomial of the form $(1 + f)$. By
Proposition \ref{prop:224}, $\varphi $ is a correct injective
homomorphism. The correspondent polynomial is constructed as
follows. Write the elements of $K$ as sums:
$$
b_i  = c_i  + d_i ,\quad c_i  \in V,\;d_i  \in \Fit(A_\Delta
),\;\,i = 1,\ldots,n.
$$
The $d_i$'s can be presented as a sum of the elements of the form:
$$
d_{ij} \frac{1}{{g_{ij} }},\quad d_{ij}  \in \Fit(A),\;g_{ij}  \in
R\smallsetminus \Delta ,\;\,i = 1,\ldots,n.
$$
Denote by $g$ the product of all the $g_{ij}$'s. Since the
polynomials $g_{ij} $'s have non-zero free terms, the polynomial
$g$ has a non-zero free term $\alpha  \ne 0$. Further, choose a
polynomial $f \in \Delta $ so that $1 + f = \frac{1}{\alpha }g$
and define the homomorphism $\varphi $ as the multiplication by
$(1 + f)$.

By the choice of the polynomial $(1 + f)$, we have $\varphi (d_i )
\in \Fit(A)$ and thus $\varphi (b_i ) \in A$, $i = 1,\ldots,n$.
Finally, since $\Fit(A_\Delta  )$ is a torsion-free module, the
homomorphism $\varphi $ is injective. \hfill $\blacksquare$

\begin{prop} \label{prop:422}
Let $A$ be a $U$-algebra and let  $A_\Delta  $ be its $\Delta
$-localisation. Then the algebras $A$ and $A_\Delta$ are
universally equivalent: \/ {\rm $\ucl(A) = \ucl(A_\Delta )$}.
\end{prop}
\proof To prove that  two algebras are universally equivalent it
suffices to prove that every finite submodel of one of the
algebras embeds into the other algebra. Clearly, every finite
submodel of the algebra  $A$ embeds into the algebra $A_\Delta$.
We, therefore, are left to prove the converse. Take a finite
submodel of the algebra $A_\Delta $: $K = \left\{ b_1 ,\ldots,b_n
\right\} $. The elements of the set  $K$ generate a finitely
generated subalgebra $B \le A_\Delta  $, which, on account of
Lemma \ref{lem:II-subofdel}, embeds into the algebra $A$.
Consequently, there exists an embedding of a finite submodel $K$
of $A$ into the algebra $A_\Delta$. \hfill $\blacksquare$

We are particularly interested in the $\Delta$-localisation of the
free metabelian Lie algebra of a finite rank $r$. Let
$c_1,\ldots,c_r $ be a linearly independent modulo $\Fit(\F)$
elements from $\F$ and let $C =  \left< c_1 ,\ldots,c_r  \right> $
be a subalgebra of $\F$. By Corollary \ref{cor:II-propofsubF}, $C$
is the free metabelian Lie algebra of the rank $r$. Although, in
general  $C \varsubsetneq \F$. However, for $\Delta$-local
algebras holds

\begin{prop} \label{prop:423}
For this notation, $C_\Delta = \F_\Delta $ and
\/ {\rm $\Fit(C_\Delta  ) = \Fit(\F_\Delta  )$}.
\end{prop}
\proof The proof is required for the second equality only. Without
loss of generality, we may assume that $c_i = a_i  + d_i ,\;i =
1,\ldots,r$, where  $a_1 ,\ldots,a_r $ is the free base of the
algebra $F$ and $d_1 ,\ldots,d_r  \in F^2 $. It suffices to show
that for an arbitrary pair of indices $(i,j)$ holds $a_i a_j \in
C_\Delta ^2 $. Let $R = k\left[x_1 ,\ldots,x_r \right]$. Write all
the products $c_i c_j$, $i > j = 1,\ldots,r$ in the form
$$
c_i c_j  = a_i a_j  \cdot f_{ij} (x_1 ,\ldots,x_r ) +
\sum\limits_{(p,q) \ne (i,j)} {a_p a_q  \cdot f_{pq} (x_1
,\ldots,x_r )} ,
$$
where $f_{ij}  \in R\smallsetminus \Delta $, $f_{pq}  \in \Delta
$. This gives a system of $C_r^2 $
 (the combination of  $r$ taken $2$ at a time) equations over $\Fit(C)$
in variables $a_i a_j ,\;\,i > j = 1,\ldots,r$. The determinant of
this system is a polynomial  $h(x_1 ,\ldots,x_r ) \in
R\smallsetminus \Delta $. Consequently, the system is compatible
over the module $\Fit_\Delta (C)$ . \hfill $\blacksquare$

\subsection{The Direct Module Extension of the Fitting's Radical of a $U$-Algebra.}

In this section we introduce another type of extension of a
$U$-algebra. As above, let  $A$ be a $U$-algebra over a field $k$.
Let $\left\{ {z_\alpha ,\;\alpha  \in \Lambda } \right\}$ be a
maximal set of linearly independent modulo $\Fit (A)$ elements
from $A$. Let $V$ denote the linear span over the field $k$ of
$\left\{ {z_\alpha ,\;\alpha  \in \Lambda } \right\}$.

Let $M$ be a torsion free module over the ring of polynomials $R$.
By the means of the module $M$ we extend the algebra $A$ to the
algebra $A \oplus M$. By the definition, the algebra $A \oplus M$
is the direct sum of $k$-vector spaces $V \oplus \Fit(A) \oplus
M$. To define the structure of an algebra on $V \oplus \Fit(A)
\oplus M$ we need to introduce the multiplication on its elements.
For the elements from $V$ the multiplication is inherited from
$A$. The multiplication of the elements from either $M$ or
$\Fit(A)$ on both elements from $\Fit(A)$ and elements from $M$
results zero. If $b \in M$ and $z_\alpha \in V$, we set $b \circ
z_\alpha = b \cdot x_\alpha$ and extend this definition of
multiplication of $b$ by elements from $V$ to be linear. This
operation turns $A \oplus M$ into a metabelian Lie algebra over
$k$.

By the definition the canonical set of generators of $A \oplus M$
is the canonical set of generators of $A$ together with an
arbitrary set of generators of $M$. One verifies that
\begin{itemize}
    \item $\Fit (A \oplus M)=\Fit(A) \oplus M$,
    \item the algebra $A \oplus M$ is a $U$-algebra.
\end{itemize}

We now point out some of the connections between the direct module
extension and the universal closure of $A$.

\begin{lem} \label{lem:431}
Let $M$ be a torsion free module over the ring of polynomials $R$.
Then for every finite tuple  $u_1 ,\ldots,u_n $ of elements from
the module $M$ and for every tuple of non-zero polynomials $f_1
,\ldots,f_n $ from $R$ there exists $u \in M$ so that  $u \cdot
f_i \ne u_i $ for any $i = 1,\ldots,n$.
\end{lem}
\proof Take a non-zero element $u_0 \in M$ and consider the
following infinite set
$$
K = \left\{ u_0 ,\;\,u_0  \cdot x_\alpha ,\;\,u_0  \cdot x_\alpha
^2 ,\;\ldots,\;u_0  \cdot x_\alpha ^m ,\;\ldots \right\},
$$
where $x_\alpha  $ is an arbitrary variable from $R$. We next show
that for any equation $u \cdot f_i = u_i$, $i=1, \ldots, n$ there
exists no more than one element from $K$ that satisfies this
equation. Let $(u_0 \cdot x_\alpha ^m ) \cdot f_i = u_i $ and $l
\ne m$. Then $(u_0 \cdot x_\alpha ^l ) \cdot f_i  = (u_0 \cdot
x_\alpha ^m ) \cdot f_i \cdot x_\alpha ^{l - m} = u_i \cdot
x_\alpha ^{l - m} \ne u_i$. The latter is implied by the  fact
that $M$ is a torsion-free $R$ module. The collection of
restrictions that are to be satisfied by $u$ is finite, while the
set $K$ is infinite. In what follows that there exists a required
element in $K$. \hfill $\blacksquare$

\begin{lem}\label{lem:432}
Let $A$ be a $U$-algebra and let $T_1 $ be a one generated
torsion-free module over the ring of polynomials $R$. Then the
algebra $A \oplus T_1 $ is $A$--discriminated by the algebra $A$.
\end{lem}
\proof Fix a finite number of non-zero elements from $A \oplus
T_1$:
$$
x_1  + u_1 ,\ldots,x_n  + u_n,\; x_i  \in A,\; u_i  \in T_1;\;u_i
= t \cdot f_i,
$$
where $t$ is the generator of  $T_1 $, $f_i  \in R$. We construct
an  $A$--homomorphism $\varphi :A \oplus T_1  \rightarrow A$ so
that $\varphi (x_i  + u_i ) \ne 0$ for any $i = 1,\ldots ,n$. For
any $a \in A$ and some $u \in \Fit(A)$, set $\varphi (a) = a$ and
$\varphi (t) = u$. This map gives rise to an $A$--homomorphism
from $A \oplus T_1 $ to $A$ in the obvious way. To show that $A
\oplus T_1 $ is $A$--discriminated by $A$ we need to choose $u \in
\Fit(A)$ so that  $x_i + u \cdot f_i \ne 0$ for any $i =
1,\ldots,n$. For indices $i$ such that $x_i  \notin \Fit(A)$ the
required inequality holds whatever $u$ is. In the event that $x_i
\in \Fit(A)$ and $f_i  = 0$ we have $\varphi (x_i ) = x_i $, i.e.
the required inequality also holds. We now use Lemma \ref{lem:431}
to choose the element $u \in \Fit(A)$ agreeing with the conditions
for the remaining indices. \hfill $\blacksquare$

\begin{lem} \label{lem:433}
For every positive integer $s$ the free module $T_s $
of the rank $s$ over the ring $R$ is discriminated by a one
generated torsion-free $R$-module $T_1 $.
\end{lem}
\proof Let $\left\{ t_1 ,\ldots,t_s \right\} $ be the free
generators of the module $T_s $. The required map has the form:
the element $t_i $ goes into $t \cdot f_i $, where $t$ is an
arbitrary non-zero element from  $T_1 $ and the $f_i $'s are
polynomials from  $R$. The polynomials $f_i $'s are chosen to take
the finite subset from  $T_s $ to non-zero elements. Consequently,
there are only finitely many restrictions imposed on the $f_i $'s.
Therefore, such polynomials exist. \hfill $\blacksquare$

\begin{prop} \label{prop:434}
If $A$ is a $U$-algebra and
 $M$ is a finitely generated module over $R$
then the algebra $A \oplus M$ is $A$--discriminated by $A$.
\end{prop}
\proof The module $M$ embeds into the free module $T_s $ of the
rank  $s$ over the ring  $R$ (see \cite{Leng}, \cite{Bourb}). This
implies that the algebra  $A \oplus M$ $A$--embeds into the
algebra $A \oplus T_s $. By Lemma \ref{lem:433} the module  $T_s $
is discriminated by the module  $T_1 $ and thus the algebra $A
\oplus T_s $ is  $A$--discriminated by the algebra $A \oplus T_1$.
According to Lemma \ref{lem:432}, the algebra $A \oplus T_1 $ is
$A$--discriminated by the algebra $A$. In what follows that the
algebra  $A \oplus M$ is $A$--discriminated by the algebra $A$.
\hfill $\blacksquare$

\begin{prop} \label{prop:435}
Let  $A$ be a $U$-algebra and let  $M$ be a finitely generated
module over $R$. Then {\rm $A- \ucl(A) = A - \ucl (A \oplus M)$}.
\end{prop}
\proof Show that every finite submodel $K = \left\{ b_1
,\ldots,b_n \right\} $ of the algebra $A \oplus M$ $A$--embeds
into the algebra $A$. Every element  $b_i  \in K$ decomposes into
the following sum: $b_i  = d_i  + u_i $, where $d_i  \in A$, $u_i
\in M$. The tuple $u_1 ,\ldots,u_n $ of elements from the module
$M$ generates a submodule  $M_0 $. Therefore, the submodel $K$
$A$--embeds into the algebra $A \oplus M_0 $, which, on account of
Proposition \ref{prop:434}, is $A$--discriminated by $A$. \hfill
$\blacksquare$

\begin{cor} \label{prop:II:ucldirext}
Let  $A$ be a $U$-algebra and let  $M$ be a finitely generated
module over $R$. Then {\rm $\ucl(A) = \ucl(A \oplus M)$}.
\end{cor}
\proof The proof is analogous to the proof of Proposition
\ref{prop:435} and left to the reader. \hfill $\blacksquare$

The $\Delta$-localisation and the direct module extension commute.

\begin{lem}
Let  $A$ be a $U$-algebra and let  $M$ be a finitely generated
module over $R$. Then $(A \oplus M)_\Delta   = A_\Delta   \oplus
M_\Delta $.
\end{lem}
\proof By the means of direct verification:
\begin{gather} \notag
\begin{split}
(A \oplus M)_\Delta   = (V \oplus \Fit(A) \oplus M)_\Delta   = V
\oplus (\Fit(A) \oplus M)_\Delta   = \\
         = V \oplus \Fit(A)_\Delta   \oplus M_\Delta   & = A_\Delta   \oplus M_\Delta .
\end{split}
\end{gather}
\hfill $\blacksquare$

For the use of algebraic geometry over the free metabelian Lie
algebra we next formulate some auxiliary statements about the
direct module extension. By $\hom _R (M, \Fit(A))$ denote the set
of all $R$-homomorphisms from $M$ to $\Fit(A)$.

\begin{lem} \label{lem:II-corresp}
In the above notation there exists a one-to-one correspondence
between the set of $R$-homomorphisms from $M$ to {\rm $\Fit(A)$}
and the set of $A$-homomorphisms from $A \oplus M$ to $A$, {\rm
$$
\hom_R (M,\Fit(A)) \leftrightarrow \hom_A (A \oplus M,A).
$$}
\end{lem}
\proof Let $\left\{ {m_\beta  } \mid \beta  \in  \Beta \right\} $
be a fixed system of generators of the module $M$. Let $\varphi
\in \hom_A (A \oplus M,A)$. In this notation, the
$A$--homomorphism $\varphi $ is completely defined by its values
on the generators $m_\beta $ of $M$. Let $\varphi (m_\beta  ) =
c_\beta  $. Since $bm_\beta   = 0$ for any $b \in \Fit(A)$, we
obtain $\varphi (bm_\beta  ) = bc_\beta   = 0$. Consequently, by
Lemma \ref{lem:214},  $c_\beta   \in \Fit(A)$. Moreover, every
module relation between the letters $m_\beta $'s rewrites as a
module relation  between the letters $c_\beta$'s. In what follows
that $\phi  =  \varphi  \mid _M $ is a module $R$-homomorphism
from $M$ to $\Fit(A)$. The converse is also true. If $\phi \in
\hom_R (M,\Fit(A))$ then there is a unique element $\varphi  \in
\hom_A (A \oplus M,A)$, corresponding to $\phi$. This
correspondence is given by the following rule:
$$
\varphi (m_\beta  ) = \phi (m_\beta  ), \; \varphi (a) = a,\; a
\in A.
$$
\hfill $\blacksquare$

\begin{lem}
Let  $M$ be a finitely generated module over $R$ then for every $m
\in M$, $m \ne 0$ there exists a homomorphism {\rm $\phi  \in
\hom_R (M,\Fit(A))$} such that $\phi (m) \ne 0$, i.e. $M$ is
approximated   by {\rm $\Fit(A)$}.
\end{lem}
\proof The module  $M$ embeds into the $R$-module $T_s $ of the
rank $s$. By Lemma \ref{lem:433}, $T_s $ is discriminated by the
module $T_1 $. At last, clearly, $T_1 $ embeds into $\Fit(A)$.
\hfill $\blacksquare$

\begin{lem} \label{lem:II-439}
Let  $M_1$ and $M_2$ be two finitely generated modules over $R$.
Then there exists a one-to-one correspondence between the set of
$R$-homomorphisms from $M_1$ to $M_2$ and the set of
$A$-homomorphisms from $A \oplus M_1$ to $A \oplus M_2$, {\rm
$$
\hom_R (M_1,M_2) \leftrightarrow \hom_A (A \oplus M_1,A \oplus
M_2).
$$}
Moreover corresponding homomorphisms {\rm $\phi \in \hom_R
(M_1,M_2)$} and {\rm $\varphi \in \hom_A (A \oplus M_1,A \oplus
M_2)$} are injective or surjective simultaneously.
\end{lem}
\proof The proof is analogous to the one of Lemma
\ref{lem:II-corresp}. \hfill $\blacksquare$

\end{document}